\documentclass[11pt]{article}
\usepackage[margin =1.0in]{geometry}
\usepackage{amssymb,amsthm}
\usepackage{amsmath}
\usepackage{enumerate}
\usepackage{hyperref}
\usepackage{cite}
\usepackage{enumitem}
\usepackage{color}
\usepackage{graphicx}
\usepackage{wrapfig}
\usepackage{url}
\usepackage{algorithmic}

\newcommand{\inclu}[0] {\ar@{^{(}->}}

\graphicspath{ {./figures/} }

\newtheorem{thm}{Theorem}[section]

\newtheorem{lem}[thm]{Lemma}
\newtheorem{cor}[thm]{Corollary}

\newtheorem{remark}{Remark}

\theoremstyle{remark}

\usepackage{mathtools}

\usepackage[boxruled]{algorithm2e}

\begin{document}
	
\title{Successive Projection for Solving Systems of Nonlinear Equations/Inequalities}
\author{Wen-Jun Zeng,\thanks{University of Michigan, Ann Arbor, MI, USA. Email: wjzeng@umich.edu}
\quad Jieping Ye\thanks{University of Michigan, Ann Arbor, MI, USA. Email: jpye@umich.edu}}	
	
\date{April 15, 2020}
\maketitle

\begin{abstract}
Solving large-scale systems of nonlinear equations/inequalities is a
fundamental problem in computing and optimization. In this paper, we propose
a generic successive projection (SP) framework for this problem.
The SP sequentially projects the current iterate onto the constraint set
corresponding to each nonlinear (in)equality. It extends
von Neumann's alternating projection for finding a point in the
intersection of two linear subspaces, Bregman's method for
finding a common point of convex sets and the Kaczmarz method
for solving systems of linear equations to the more general case of multiple nonlinear
and nonconvex sets. The existing convergence analyses on randomized Kaczmarz
are merely applicable to linear case. There are no theoretical convergence results of the SP
for solving nonlinear equations. This paper presents the
first proof that the SP locally converges to a solution of nonlinear
equations/inequalities at a linear rate. Our work establishes
the convergence theory of the SP for the case of multiple nonlinear and nonconvex sets.
Besides cyclic and randomized projections, we devise two new
greedy projection approaches that significantly
accelerate the convergence. Furthermore, the theoretical
bounds of the convergence rates are derived. We reveal that the convergence
rates are related to the Hoffman constants of the Jacobian matrix
of the nonlinear functions at the solution. Applying the SP to solve the graph realization problem,
which attracts much attention in theoretical computer science, is discussed.
\end{abstract}

\section{Introduction}

Solving systems of nonlinear equations is a foundational problem in computing, numerical analysis
and optimization. Nonlinear equations are ubiquitous since many physical phenomena are essentially
nonlinear. Solving nonlinear equations arises from numerically solving partial
differential equations (PDEs) and unconstrained optimization problem \cite{Freund}.
In unconstrained multivariate optimization, finding a minimizer of an objective function
can be converted into solving the system of nonlinear equations in which the gradient equals the
zero vector. Developing efficient solvers both in practice and theory for large-scale
nonlinear equations is an important topic of computer science because
it has many applications in engineering and science. A noticeable example highly
related to solving nonlinear equations is the graph realization problem (GRP) \cite{So:SODA}, which
has received a great deal of attention in theoretical computer science. The GRP
can be formulated as follows. Given a subset of pairwise Euclidean distances, find
the coordinates of points that match those distances \cite{So:SODA,Singer}. The GRP amounts to solving
a system of quadratic equations. It has a variety of applications such as surveying,
satellite ranging, sensor network localization, molecular conformation and structural
biology. Another example of quadratic equations is the phase retrieval, which aims at
recovering a complex-valued signal from magnitude-only measurements \cite{WF}. It has a wide
range of applications in cryo-electron microscopy and imaging science \cite{Nature}.
Solving nonlinear equations with some specific forms, e.g., systems of polynomial
equations \cite{Poly:Eq}, also attracts the attention of the theoretical computer
scientists.

On the other hand, solving systems of nonlinear inequalities or combinations of nonlinear equalities and inequalities
is central to constrained optimization. In nonlinear optimization, not only equality constraints
but also inequality ones can be encountered. Finding a feasible point or determining the feasible
region of a constrained optimization problem, which is referred to as \emph{feasibility problem},
is equivalent to solving a combined system of nonlinear equalities and inequalities. In addition,
solving the Karush--Kuhn--Tucker (KKT) conditions \cite{Bertsekas}, i.e., the necessary optimality condition of
constrained optimization problems, also amounts to solving this problem. Compressed sensing \cite{Candes:SPM} and matrix
completion are two research topics that attract much attention in theoretical computer science,
data science, signal processing and statistics in recent years \cite{Hardt,Ruoyu,Jiang}. The main constraints in compressed
sensing and matrix completion are the sparsity constraint and low-rank constraint, respectively, which
are nonlinear and nonconvex. The two problems in fact can be formulated as feasibility problems.
For example, the reader is referred to \cite{Jiang} on how to cast the matrix completion as a
feasibility problem.

\subsection{Prior Work}

Numerous methods have been developed for solving nonlinear equations.
When nonlinear equations reduces to linear ones, there are
many efficient solvers for this simpler problem, including the Jacobi,
Gauss--Seidel and the successive over-relaxation (SOR) methods \cite{Greenbaum,Golub}. When
the coefficient matrix is symmetric and positive-definite, the conjugate
gradient method \cite{CG}, which belongs to the Krylov subspace methods
\cite{Greenbaum,Golub,Trefethen}, is an efficient solver. The best known time complexity
for solving a general linear system is $\mathcal{O}(n^{2.373})$, where $n$ is
the number of unknowns variables \cite{VVW}. Besides the general linear equations,
faster algorithms for solving some special forms of linear systems, e.g.,
the symmetric diagonally dominant (SDD) linear systems, are of great interest
in theoretical computer science. Several nearly-linear time solvers have been
developed recently for the SDD systems \cite{Spielman,Koutis,Kelner,YTLee}.

A remarkable iterative algorithm suitable for solving larger-scale linear systems
is the Kaczmarz method \cite{Kaczmarz}, which was first proposed by Kaczmarz and
later re-discovered in the field of image reconstruction with the name of algebraic
reconstruction technique (ART). In practice, the Kaczmarz method can outperform the
conjugate gradient method. Therefore, it is widely used in practical applications, such as
computerized tomography (CT) due to its high efficiency. The idea behind the Kaczmarz method
is the \emph{alternating projection (AP)} or \emph{successive projection (SP)}\footnote{It is called
alternating projection when there are two constraints while successive projection is used
for multiple constraints.} \cite{Neumann,Bregman,AP:Book,Bauschke}. Alternating projection was applied to
find a point in the intersection of two linear subspaces and von Neumann proved the convergence under
this simple setting \cite{Neumann}. Geometrically, each linear equation of the system is a hyperplane.
At each iteration, the Kaczmarz method picks one hyperplane and then projects the current solution onto the this hyperplane.
The original Kaczmarz method selects the hyperplane cyclically. The performance of the cyclic Kaczmarz
may become worse if there is an undesirable ordering of the rows of the coefficient matrix.
Instead of the cyclic order, Strohmer and Vershynin proposed a randomized Kaczmarz method
in a seminal work \cite{Strohmer}, which picks a hyperplane randomly at each iteration.
The randomized Kaczmarz rule makes the convergence analysis easier compared to its
cyclic counterpart. The coordinate descent method \cite{Nesterov,Wright} in optimization is highly
related to the Kaczmarz method. Based on such a relation, an improved asymptotic convergence
guarantee for the Kaczmarz method is obtained via the accelerated coordinate descent method \cite{YTLee}.

Solving nonlinear equations is much more difficult than solving linear ones.
The most classical solver for the nonlinear case is the Newton's method.
At each iteration of the Newton's method, it requires to solve a linear system,
where the coefficient matrix is the Jacobian matrix at the current solution, to obtain the update
direction, i.e., the Newton direction. In other words, the Newton's method
converts a nonlinear problem into a series of linear ones. However, computing
the Jacobian can be computationally demanding, especially for large-scale problems.
The high per-iteration complexity constitutes one main drawback of the Newton's method.
Moreover, when the Jacobian is singular or ill-conditioned, convergence of the Newton's
method is questionable. Recently, an accelerated residual method for solving nonlinear
equations \cite{Freund} has been developed by exploiting Nesterov's accelerated gradient method
\cite{Nesterov2} for convex optimization.

Compared to solving nonlinear equations, solving nonlinear inequalities can be more challenging.
The AP or SP is a powerful tool suited for this problem. Bregman extended von Neumann's and
Kaczmarz's methods that are limited to linear subspces/equations to find a common point of convex sets.
The existing proofs of the convergence of the AP or SP \cite{Neumann,Bregman,AP:Book,Bauschke},
including the Kaczmarz method \cite{Strohmer,Kelner,YTLee}, explicitly or implicitly exploit the non-expansiveness
of the projection onto convex sets (POCS). As a result, these convergence analyses are not applicable
to the case of nonlinear equations, where the constraint sets associated with nonlinear equations
are nonconvex in general. Very few theoretical results on the convergence of the SP are available
except the one based the metric regularity of two manifolds \cite{Lewis}. However, the convergence
result of \cite{Lewis} is limited to two sets and cannot handle multiple nonconvex sets.
Different from \cite{Lewis}, the analysis techniques developed in our work are applicable to
multiple nonconvex sets.

\subsection{Summary of Contributions}

Aiming to overcome the high per-iteration complexity and other drawbacks of the Newton's method,
this paper proposes an SP framework for solving systems of nonlinear equations/inequalities.
The SP is simple and easy to implement: it just sequentially projects the current iterate
onto the constraint set associated with each nonlinear (in)equality.

We summarize the contributions of this paper as follows.
\begin{itemize}
  \item A generic SP framework including seven variants, namely, cyclic, uniformly and
  non-uniformly random, randomly permuted, greedy and normalized greedy, and mean projections,
  is proposed to solve systems of nonlinear equations/inequalities.
  \item The SP generalizes von Neumann's alternating projection for finding a point
  in the intersection of two linear subspaces, Bregman's successive projection
  for finding a common point of convex sets, and the Kaczmarz method for solving systems
  of linear equations.
  \item We establish the theory of convergence and iteration complexity of the SP method for multiple nonlinear/nonconvex sets.
  We give the first proof that the SP locally converges to a solution of nonlinear equations/inequalities at a linear rate.
  Most previous theoretical convergence analyses on alternating projection are limited to linear and (or) convex cases.
  \item Two accelerated versions of the SP using greedy selection rule are devised.
  To the best of our knowledge, the two greedy variants have the fastest convergence
  rate among all projection based approaches. We give tighter bounds on the
  convergence rates of the greedy projection methods. That is, we prove that the
  two greedy methods are faster than other variants in theory.
  \item In numerous applications in engineering and science where solving nonlinear systems is required,
  such as graph realization, sensor network localization, matrix completion,
  compressed sensing, phase retrieval and molecular conformation, the projections
  onto the nonlinear constraint sets have closed-form solutions or are computationally
  tractable. In these cases, the SP is efficient and highly competitive compared
  to other methods.
\end{itemize}

In addition, our method can be applied to compute the projection onto the constraint sets
in projected gradient method, especially for nonconvex optimization problems.
We hope that our work will provide new insights to solving nonlinear systems in both theory and practice
and achieve successes in practical applications.

\subsection{Paper Organization and Notations}

The remainder of this paper is organized as follows. The problem is
formulated and the preliminaries are introduced in Section \ref{Sec:Preliminaries}.
In Section \ref{Sec:SP}, we present the generic SP framework for solving systems
of nonlinear equations/inequalities and detail its variants. The main results of
the convergence of the SP are presented in Section \ref{Sec:MainResult}.
Section \ref{Sec:Convergence} gives detailed proofs of these theoretical results.
Numerical results are provided in Section \ref{Sec:Simulation}
to demonstrate the fast convergence rate and efficiency of the proposed methods.
Some stuff are in the appendices.

Bold capital upper case and lower case letters represent matrices
and vectors, respectively. The identity matrix is written as $\pmb I$.
The superscripts $(\cdot)^\top$, $(\cdot)^\mathrm{H}$ and $(\cdot)^\dagger$
denote the transpose, Hermitian transpose and Moore-Penrose pseudoinverse, respectively.
$\mathbb{E}[\cdot]$ is expectation. The $\ell_p$-norm of a vector is represented by $\|\cdot\|_p$
while the spectral norm and Frobenius norm of a matrix are denoted as
$\|\cdot\|_2$ and $\|\cdot\|_\mathrm{F}$, respectively. Finally, $\mathbb{R}$, $\mathbb{R}_+$
and $\mathbb{C}$ stand for the sets of real, non-negative real, and complex numbers, respectively.

\section{Problem Formulation}\label{Sec:Preliminaries}

We consider solving the system of nonlinear equations
\begin{equation}\label{NLE}
    \pmb f(\pmb x) = \pmb 0
\end{equation}
where $\pmb x =[x_1,\cdots, x_n]^\top\in\mathbb{R}^n$ collects the unknown variables
and $\pmb f:\mathbb{R}^n \rightarrow\mathbb{R}^m$ is a vector-valued function, which
can be written as $\pmb f(\pmb x)=[f_1(\pmb x),\cdots, f_m(\pmb x)]^\top$ with
$f_i:\mathbb{R}^n \rightarrow\mathbb{R}$ ($i=1,\cdots,m$) being a real-valued function.
In this paper, we focus on the case that the number of equations $m$ is not less than
the number of unknowns $n$, i.e., $m\ge n$. The proposed SP method is also applicable
to the under-determined case with $m<n$. We also consider solving the system of nonlinear
inequalities
\begin{equation}\label{NLInE}
    \pmb f(\pmb x) \le \pmb 0
\end{equation}
or combinations of nonlinear equalities and inequalities. Solving \eqref{NLInE} is
equivalent to finding a feasible point of the nonlinear constraints $\{\pmb x | \pmb f(\pmb x) \le \pmb 0\}$,
which is a central problem in nonlinear optimization. When $\pmb f$ takes the affine form
$\pmb f(\pmb x) = \pmb A\pmb x - \pmb b$, \eqref{NLE} reduces to a system of linear equations
$\pmb A\pmb x = \pmb b$ and \eqref{NLInE} becomes a linear feasibility problem $\pmb A\pmb x \le \pmb b$
that lies central in linear programming. Some nonlinear equations/inequalites frequently
encountered in engineering and science are listed in Appendix A.

We detail the GRP as another application example of solving nonlinear equations.
Given the coordinates of $n_v$ points in $\mathbb{R}^d$ with $d\ge 1$
being the dimension of the coordinate space, computing the distance
between any two points is an easy task. Graph realization is the inverse problem: given a subset of pairwise
distances, find the coordinates of points in $\mathbb{R}^d$
that match those distances is never easy. This problem has been proved to be NP-hard
for any $d\ge 1$. The GRP is formulated as follows.
Given a graph $G = (V;E)$ with $n_v$ vertices and
$m$ edges, the pairwise Euclidean distances $\{r_{ij}\}_{(i,j)\in E}$, decide if there exist
vectors $\pmb x_1,\cdots,\pmb x_{n_v}\in \mathbb{R}^d$ such that
\begin{equation}\label{GRP}
    \|\pmb x_i-\pmb x_j\|_2=r_{ij},~\mathrm{for~all}~{(i,j)\in E}.
\end{equation}
The configuration $\pmb x = [\pmb x_1^\top,\cdots,\pmb x_{n_v}^\top]^\top\in \mathbb{R}^n$
with $n=n_vd$ is called a realization of $G$. It is clear that the
graph realization problem is equivalent to judging whether the quadratic
equations of \eqref{GRP} has a solution and then solving it.

\section{Successive Projection Methods}\label{Sec:SP}

Denote the constraint set associated with the $i$th nonlinear equation as
\begin{equation}\label{Si:def}
    \mathcal{S}_i = \{\pmb x|f_i(\pmb x)=0\}, \quad i = 1,\cdots,m
\end{equation}
then the solution set of \eqref{NLE}:
\begin{equation}\label{Set:NLE}
    \mathcal{S}=\{\pmb x|\pmb f(\pmb x)=\pmb 0\}
\end{equation}
can be expressed as
\begin{equation}\label{S:Si}
    \mathcal{S}=\cap_{i=1}^m \mathcal{S}_i.
\end{equation}
It is clear that solving the nonlinear equations of \eqref{NLE} is equivalent
to solving the following feasibility problem
\begin{equation}\label{S:Si}
    \mathrm{find}~\pmb x \in \cap_{i=1}^m \mathcal{S}_i.
\end{equation}
The SP is an iterative algorithm for finding a feasible point of $\mathcal{S}=\cap_{i=1}^m \mathcal{S}_i$.
Denote the result of the $k$th iteration as $\pmb x^k$. At the $k$th iteration, the SP projects the current iterate $\pmb x^k$
onto the $i_k$th ($i_k\in\{1,\cdots,m\}$) constraint set $\mathcal{S}_{i_k} =\{\pmb x|f_{i_k}(\pmb x)=0\}$.
Here the projection of an arbitrary point $\pmb z\in \mathbb{R}^n$ onto the set $\mathcal{S}_i$ is
the point in $\mathcal{S}_i$ that is ``closest'' (in Euclidean distance) to $\pmb z$, which is defined as
\begin{equation}\label{Proj:Si}
    \Pi_i(\pmb z) := \arg\min_{\pmb x \in \mathcal{S}_i} \|\pmb x - \pmb z\|_2^2.
\end{equation}
The iteration of the SP is then expressed as
\begin{equation}\label{}
    \pmb x^{k+1} = \Pi_{i_k}(\pmb x^k).
\end{equation}
The SP is listed in Algorithm \ref{Algo:SP}. When there are two sets and the
sets are linear subspace, the SP reduces to the von Neumann's
method of alternating projection. When the equations are linear,
the SP reduces to the Kaczmarz method.
\begin{algorithm}
    \caption{SP for Solving Nonlinear Equations}\label{Algo:SP}
    \algsetup{indent=2em}
    \begin{algorithmic}
    \vspace{1ex}
    \STATE \textbf{Initialization:} Choose $\pmb x^0\in\mathbb{R}^n$.
    \FOR{$k=0,1,\cdots,$}
    \STATE Choose index $i_k\in\{1,\cdots,m\}$
    \STATE $\pmb x^{k+1} = \Pi_{i_k}(\pmb x^k)$
    \STATE \textbf{Stop} if convergence condition is satisfied.
    \ENDFOR
    \end{algorithmic}
\end{algorithm}

Several rules to select $i_k$ are considered as follows, which result in different projection methods.
\begin{itemize}
  \item Cyclic Projection (CP): $i_k$ cyclically takes value from $\{1,\cdots,m\}$.
  Every $m$ iterations are called one cycle. The original Karczmarz method for
  solving linear systems of equations \cite{Kaczmarz} adopts the
  cyclic update rule. This rule is also used in the Gauss-Seidel
  method for solving linear equations \cite{Greenbaum} and the coordinate descent method
  for unconstrained minimization problems \cite{Wright}, where each coordinate is updated with a cyclic order.
  \item Random Projection (RP): $i_k$ is randomly chosen from $\{1,\cdots,m\}$ with equal probability.
  \item Randomly Permuted Projection (RPP): like the CP, the RPP iterates cycle-by-cycle.
  But at each cycle of the RPP, the order of projections is a random permutation of $\{1,\cdots,m\}$.
  The set of the selected indices of one cycle of the RPP is denoted as $\{i_1,\cdots,i_m\}$, which
  is a random permutation of $\{1,\cdots,m\}$.
  \item Non-uniformly Random Projection (NRP): if the $\ell_2$-norms of the gradients
  of each nonlinear equation at the solution $\pmb x^*$
are available, $i_k$ is non-uniformly sampled from $\{1,\cdots,m\}$
  with probability of
  \begin{equation}\label{NRP:rule}
  \frac{\|\nabla f_i(\pmb x^*)\|_2^2}{\sum_{i=1}^m\|\nabla f_i(\pmb x^*)\|_2^2},
  \quad i=1,\cdots,m
  \end{equation}
  where $\nabla f_i(\pmb x^*)$ is the gradient of $f_i(\cdot)$ at the solution $\pmb x^*$.
  Since the solution $\pmb x^*$ is unknown, we cannot evaluate
$\|\nabla f_i(\pmb x^*)\|_2$ in general. However, $\|\nabla f_i(\pmb x^*)\|_2$
can be computed without knowing $\pmb x^*$ in many problems. Some examples are listed in Appendix B.
We see that the RP treats each equation equally while the NRP evaluates the importance of
each equation based on their gradients.

\item Greedy Projection (GP): we propose to select $i_k$ by
\begin{equation}\label{GP:rule}
    i_k= \arg\max_{1\le i\le m} \left|f_i(\pmb x^k)\right|.
\end{equation}
Clearly, the greedy projection chooses $i_k$ with the maximum magnitude of the residual.
Therefore, computing the residual is required at each iteration while
it does not need for the random and cyclic projections. We will provide both theoretical
analysis and experimental results to validate that the GP with the maximum residual
rule converges faster than the CP and RP at the expense of an extra residual evaluation.
\item Normalized Greedy Projection (NGP): $i_k$ is selected by
\begin{equation}\label{NGP:rule}
    i_k= \arg\max_{1\le i\le m}
    \frac{\left|f_i(\pmb x^k)\right|}{\|\nabla f_i(\pmb x^*)\|_2}.
\end{equation}
Compared to the GP, the NGP incorporates the $\ell_2$-norms of the gradients.
If the magnitudes of $\{\nabla f_i(\pmb x)\}_i$ are significantly
different with each other, the GP may not be optimal. By normalizing
the residual with the $\ell_2$-norm of the gradient, the NGP eliminates
the magnitude difference and can perform better than the GP. When
$\|\nabla f_i(\pmb x^*)\|_2$ does not depend on $i$, just as the unsquared
circle and elliptic equations, the NPG is equivalent to the GP.

\item Mean Projection (MP): at each iteration, the MP computes all projections
onto $m$ constraint sets and then average them, yielding the following iteration
\begin{equation}\label{}
    \pmb x^{k+1} = \frac{1}{m}\sum_{i=1}^m\Pi_i(\pmb x^k).
\end{equation}
It is obvious that the MP requires to compute $m$ projections in one iteration
while the CP, RP and GP only need one.
\end{itemize}
If the projection onto each nonlinear constraint set is computationally tractable,
the SP is efficient. Some examples where the projection has closed-form expressions
or can be easily computed are listed in Appendix C. These examples
are also frequently encountered in science and engineering applications.
For convex sets, it is known that the POCS is unique and non-expansive
\cite{Bertsekas,Bauschke}. Most existing theoretical convergence analyses of the SP
are based on the uniqueness and non-expansiveness of POCS, which makes them merely applicable
to convex sets. The convergence analysis of the SP for nonconvex sets is much more difficult
without the two properties of POCS.

\section{Main Results}\label{Sec:MainResult}

In the convergence analysis, we only need to assume that the projections $\{\Pi_i(\cdot)\}_{i=1}^m$
are differentiable at a solution $\pmb x^*$ with $\pmb f(\pmb x^*)=\pmb 0$ or equivalently
$\pmb x^*\in \cap_{i=1}^m \mathcal{S}_i$. This assumption is quite mild and often holds true
in many practical applications. For example, it is easy to check that the four projections
of \eqref{Proj:A}--\eqref{Proj:Ei} are differentiable at any solution $\pmb x^*$.

To facilitate to describe our main results and the basic technique
used in the convergence proof, we introduce a mapping $\Phi$:
$\mathbb{R}^n\rightarrow \mathbb{R}^n$, which describes the operation of
one basic iteration of variants of SP. For MP, $\Phi(\cdot)$ has the form
\begin{equation}\label{Phi:MP}
    \Phi_\mathrm{MP}(\cdot) = \frac{1}{m}\sum_{i=1}^m\Pi_i(\cdot)
\end{equation}
while for CP, $\Phi(\cdot)$ describes a cycle of $m$ projections, which is
\begin{equation}\label{Phi:CP}
    \Phi_\mathrm{CP}(\cdot) = \Pi_m\Pi_{m-1}\cdots\Pi_2\Pi_1(\cdot).
\end{equation}
Similar to the CP, $\Phi(\cdot)$ of RPP is a cycle of randomly permuted projections:
\begin{equation}\label{Phi:RPP}
    \Phi_\mathrm{RPP}(\cdot) = \Pi_{i_m}\Pi_{i_{m-1}}\cdots\Pi_{i_2}\Pi_{i_1}(\cdot)
\end{equation}
where $\{i_1,\cdots,i_m\}$ is a random permutation of $\{1,\cdots,m\}$.
For the RP and NRP, $\Phi_\mathrm{RP}(\cdot) = \Pi_{i_k}(\cdot)$ with $i_k$ being
randomly selected from $\{1,\cdots,m\}$. $\Phi_\mathrm{GP}(\cdot)$
and $\Phi_\mathrm{NGP}(\cdot)$ have the same form as $\Phi_\mathrm{RP}(\cdot)$ but $i_k$ is chosen according
to \eqref{GP:rule} and \eqref{NGP:rule}, respectively. Sometimes we omit the subscript and just use $\Phi(\cdot)$.
The following unified form describes the operation of one iteration\footnote{It
is one iteration for MP, RP, NRP, GP and NGP while it is a cycle ($m$ iterations) for CP and RPP.}
of the SP
\begin{equation}\label{IterFormat:SP}
\pmb x^{k+1} = \Phi(\pmb x^k),~k=0,1,\cdots.
\end{equation}
Note that $\pmb x^*$ is a fixed-point of the mapping $\Phi(\cdot)$
when $\Phi(\cdot)$ take forms of $\Phi_\mathrm{MP}(\cdot)$, $\Phi_\mathrm{CP}(\cdot)$,
$\Phi_\mathrm{RPP}(\cdot)$, $\Phi_\mathrm{RP}(\cdot)$, $\Phi_\mathrm{NRP}(\cdot)$,
$\Phi_\mathrm{GP}(\cdot)$ and $\Phi_\mathrm{NGP}(\cdot)$, i.e.,
\begin{equation}\label{x:star:FP}
\pmb x^* = \Phi(\pmb x^*)
\end{equation}
because of $\pmb x^* = \Pi_i(\pmb x^*)$ for $i=1,\cdots,m$, which is
due to $\pmb x^* \in \mathcal{S}_i$. As a result, the
SP of \eqref{IterFormat:SP} can be viewed as a fixed-point iteration.

We are ready to formally state the convergence results.
\begin{thm}[Convergence of MP, CP and RPP]\label{Thm:Converge:CP}
The variants of the SP of Algorithm \ref{Algo:SP} with cyclic and randomly permuted
index selection rules, i.e., the CP and RPP, as well as the MP
locally converge to a solution of the system of nonlinear equations
$\pmb f(\pmb x)=\pmb 0$, which is denoted as $\pmb x^*$, at a linear
rate if $\mathrm{rank}(\pmb U)=n$ with $\pmb U$ being defined as
\begin{equation}\label{U:def}
     \pmb U = \left[\frac{\nabla f_1(\pmb x^*)}{\|\nabla f_1(\pmb x^*)\|_2},\cdots,
     \frac{\nabla f_m(\pmb x^*)}{\|\nabla f_m(\pmb x^*)\|_2}\right]\in\mathbb{R}^{n\times m}.
\end{equation}
Let $\{\pmb x^k\}_{k=1,2,\cdots}$ be the sequence generated by Algorithm \ref{Algo:SP}.
There exists a neighborhood centered at $\pmb x^*$ with radius $\delta>0$
\begin{equation}\label{Neighborhood}
    \mathcal{B}_\delta(\pmb x^*) =
    \left\{\pmb x\left|\|\pmb x - \pmb x^*\|_2 < \delta\right\} \right.
\end{equation}
such that if the initial guess $\pmb x^0\in \mathcal{B}_\delta(\pmb x^*)$, then
\begin{equation}\label{}
    \|\pmb x^k - \pmb x^*\|_2 < \gamma^k\|\pmb x^0 - \pmb x^*\|_2
\end{equation}
where the convergence rate
\begin{equation}\label{}
    \gamma = \epsilon + \|\nabla \Phi(\pmb x^*)\|_2
\end{equation}
satisfies $0<\gamma<1$ since the spectral norm $\|\nabla \Phi(\pmb x^*)\|_2<1$
and $\epsilon>0$ can be arbitrarily small as $\delta$ decreases.
\end{thm}

As the iteration progresses, it has $\mathop{\lim}\limits_{k\rightarrow\infty}\epsilon=0$
due to $\mathop{\lim}\limits_{k\rightarrow\infty}\pmb x^k=\pmb x^*$
and $\mathop{\lim}\limits_{k\rightarrow\infty}\delta=0$. Since $\epsilon$ is small
enough, $\|\nabla \Phi(\pmb x^*)\|_2$ is dominant in the convergence rate $\gamma$.
We focus on the asymptotic convergence rate
\begin{equation}\label{}
    \mathop{\lim}\limits_{k\rightarrow\infty}\gamma
    = \mathop{\lim}\limits_{k\rightarrow\infty}(\epsilon + \|\nabla \Phi(\pmb x^*)\|_2)
    =\|\nabla \Phi(\pmb x^*)\|_2.
\end{equation}

\begin{cor}[Asymptotic Convergence Rate of MP]\label{Thm:Rate:MP}
The asymptotic convergence rate of the MP is
\begin{equation}\label{}
    \mathop{\lim}\limits_{k\rightarrow\infty}\gamma_\mathrm{MP}
    =\|\nabla\Phi_\mathrm{MP}(\pmb x^*)\|_2
    =\sqrt{1 - \frac{1}{m}\sigma_{\min}^2(\pmb U)}
\end{equation}
where $\sigma_{\min}(\pmb U)$ is the minimum singular value of $\pmb U$.
\end{cor}

Noting that $\|\pmb U\|_\mathrm{F}^2=m$ since the $\ell_2$-norm of each column of $\pmb U$ is one,
the asymptotic convergence rate of the MP can be expressed in terms of the condition number
\begin{equation}\label{AsymRate:MP}
    \mathop{\lim}\limits_{k\rightarrow\infty}\gamma_\mathrm{MP}
    =\sqrt{1 - \frac{1}{\kappa^2(\pmb U)}}
\end{equation}
where the condition number is defined as
\begin{equation}\label{kappa:def}
   \kappa(\pmb U)= \frac{\|\pmb U\|_\mathrm{F}}{\sigma_{\min}(\pmb U)}.
\end{equation}

Since the RP and NRP are randomized algorithms, we present their convergence in expectation.

\begin{thm}[Expectation Convergence of RP and NRP and Asymptotic Convergence Rates]\label{Thm:Rate:RP}
When $i_k$ is uniformly sampled from $\{1,\cdots,m\}$ for RP and non-uniformly sampled
according to the distribution of \eqref{NRP:rule} for NRP, the RP and NRP locally converge
to a solution $\pmb x^*$ at a linear rate in expectation if
$\mathrm{rank}(\pmb U)=n$ and $\pmb x^0\in \mathcal{B}_\delta(\pmb x^*)$ with $\mathcal{B}_\delta(\pmb x^*)$
defined in \eqref{Neighborhood}. In particular, let $\{\pmb x^k\}_{k=0,1,\cdots}$ be the sequence
generated by the RP or NRP, we have
\begin{equation}\label{}
   \mathbb{E}_{i_k}\left[\|\pmb x^{k+1} - \pmb x^*\|_2\right]
    < \gamma_\mathrm{RP}^{k+1}\|\pmb x^0 - \pmb x^*\|_2
\end{equation}
for the RP and
\begin{equation}\label{}
   \mathbb{E}_{i_k}\left[\|\pmb x^{k+1} - \pmb x^*\|_2\right]
    < \gamma_\mathrm{NRP}^{k+1}\|\pmb x^0 - \pmb x^*\|_2
\end{equation}
for the NRP, where $\gamma_\mathrm{RP},\gamma_\mathrm{NRP}\in(0,1)$ are
the convergence rates of the RP and NRP, respectively.
Furthermore, the asymptotic convergence rates are
\begin{equation}\label{AsymRate:RP}
    \mathop{\lim}\limits_{k\rightarrow\infty}\gamma_\mathrm{RP}
    =\sqrt{1 - \frac{1}{\kappa^2(\pmb U)}}
\end{equation}
and
\begin{equation}\label{AsymRate:NRP}
    \mathop{\lim}\limits_{k\rightarrow\infty}\gamma_\mathrm{NRP}
    =\sqrt{1 - \frac{1}{\kappa^2(\pmb G)}}
\end{equation}
where
\begin{equation}\label{}
    \pmb G=[\nabla f_1(\pmb x^*),\cdots,\nabla f_m(\pmb x^*)]\in \mathbb{R}^{n\times m}
\end{equation}
is the transpose of the Jacobian matrix of $\pmb f(\cdot)$ at the solution
$\pmb x^*$ and $\kappa(\cdot)$ is the condition number with definition in \eqref{kappa:def}.
\end{thm}

\begin{thm}[Convergence of NGP and GP and Convergence Rates]\label{Thm:Rate:GP}
The GP and NGP with index selection rules in \eqref{GP:rule} and \eqref{NGP:rule}, respectively,
locally converge to a solution $\pmb x^*$ at a linear rate if $\mathrm{rank}(\pmb U)=n$ and
$\|\pmb x^0 - \pmb x^*\|_2 < \min(\delta,\delta')$, where $\delta'$ is another radius parameter
making Lemma \ref{Lemma:Greedy:Gradient} hold true. Specifically, let $\{\pmb x^k\}_{k=1,2,\cdots}$
be the sequence generated by the NGP or GP, we have
\begin{equation}\label{}
   \|\pmb x^k - \pmb x^*\|_2
    < \gamma_\mathrm{NGP}^k\|\pmb x^0 - \pmb x^*\|_2
\end{equation}
for NGP and
\begin{equation}\label{}
   \|\pmb x^k - \pmb x^*\|_2
    <\gamma_\mathrm{GP}^k\|\pmb x^0 - \pmb x^*\|_2
\end{equation}
for GP, where $\gamma_\mathrm{NGP},\gamma_\mathrm{GP}\in(0,1)$ are
the convergence rates of NGP and GP, respectively. Moreover, their
asymptotic convergence rates are
\begin{equation}\label{AsymRate:NGP}
    \mathop{\lim}\limits_{k\rightarrow\infty}\gamma_\mathrm{NGP}
    =\sqrt{1-h_\infty^2(\pmb U^\top)}
\end{equation}
and
\begin{equation}\label{}
    \mathop{\lim}\limits_{k\rightarrow\infty}\gamma_\mathrm{GP}
    =\sqrt{1-\frac{h_\infty^2(\pmb G^\top)}{\left\|\pmb G^\top\right\|_{2,\infty}^2}}
\end{equation}
where $h_\infty(\pmb U^\top)$ is the Hoffman type
constant \cite{Hoffman} of $\pmb U^\top$, which is defined as
\begin{equation}\label{}
  h_p(\pmb U^\top) = \inf_{\pmb v\ne \pmb 0}\frac{\|\pmb U^\top \pmb v\|_p}{\|\pmb v\|_2}
\end{equation}
with $\|\cdot\|_p$ denoting the $\ell_p$-norm ($1\le p\le\infty$), and
$\left\|\pmb G^\top\right\|_{2,\infty}$ is the $\ell_{2,\infty}$-norm of
$\pmb G^\top$ with definition being
\begin{equation}\label{}
\left\|\pmb G^\top\right\|_{2,\infty}=\max_{1\le i\le m}\|\nabla f_i(\pmb x^*)\|_2.
\end{equation}
\end{thm}

As a note, the NGP is faster than the RP due to
$h_\infty^2(\pmb U^\top)\ge\frac{1}{m}\sigma_{\min}^2(\pmb U^\top)$.
Note that the worst case where the equality $h_\infty^2(\pmb U^\top)
=\frac{1}{m}\sigma_{\min}^2(\pmb U^\top)$ holds seldom occurs\footnote{Detailed
explanation on this point can be found in the proof of Theorem \ref{Thm:Rate:GP}.}.
Numerical experiment results in Section \ref{Sec:Simulation}
show that the NGP is much faster than the RP in general.

\section{Proofs of Convergence}\label{Sec:Convergence}

To prove the convergence results in Section \ref{Sec:MainResult}, we need the following preliminaries.

\subsection{Differential of the Mapping and Spectral Norm}

The concept of tangent space \cite{Jeffrey} is useful in the convergence analysis.
The set $\mathcal{S}_i = \{\pmb x \left|f_i(\pmb x)=0\right.\}$ defines a manifold
on $\mathbb{R}^n$. The tangent space of $\mathcal{S}_i$ at $\pmb x\in\mathcal{S}_i$
is defined as
\begin{equation}\label{TanSpace:def}
    T_{\mathcal{S}_i}(\pmb x)
    =\left\{\pmb v\in \mathbb{R}^n \left|\nabla f_i(\pmb x)^\top\pmb v=0\right.\right\}
\end{equation}
where $\nabla f_i(\pmb x)$ is the gradient of $f_i(\pmb x)$. We see
that any vector in the tangent space is orthogonal to the gradient.

\begin{lem}[Jacobian Matrix of Projection \cite{Lewis,Robinson}]\label{Jacob:Proj}
The Jacobian matrix of the projection operator $\Pi_i(\cdot)$ of \eqref{Proj:Si}
at a point $\pmb x\in\mathcal{S}_i$ equals the the projection onto the tangent
space of $\mathcal{S}_i$ at $\pmb x$, i.e.,
\begin{equation}\label{JPi:Tan}
\nabla\Pi_i(\pmb x) = \Pi_{T_{\mathcal{S}_i}(\pmb x)}.
\end{equation}
\end{lem}

\emph{Proof.} See a proof in \cite{Lewis} where the main technique
is based on the derivatives of the mapping functions developed by Robinson \cite{Robinson}.
\hfill $\square$

Since the tangent space is a subspace in $\mathbb{R}^n$, the
projection $\Pi_{T_{\mathcal{S}_i}(\pmb x)}$ is linear and
can be represented as projection matrices. By the definition of
the tangent space of \eqref{TanSpace:def}, we know that $T_{\mathcal{S}_i}$ is
the kernel (null) space of $\nabla f_i(\pmb x)^\top$, i.e.,
$T_{\mathcal{S}_i}=\mathrm{ker}(\nabla f_i(\pmb x)^\top)$.
As a result, $\Pi_{T_{\mathcal{S}_i}(\pmb x)}$ is the projection
matrix onto the orthogonal complement of the range space spanned by
$\nabla f_i(\pmb x)$ and can be computed analytically as
\begin{equation}\label{Proj:Tan:Si}
    \Pi_{T_{\mathcal{S}_i}(\pmb x)}
    = \pmb I - \frac{\nabla f_i(\pmb x)\nabla f_i(\pmb x)^\top}{\|\nabla f_i(\pmb x)\|_2^2}
\end{equation}
where $\pmb I$ is the identity matrix. Herein, it assumes $\nabla f_i(\pmb x)\ne\pmb 0$.
If $\nabla f_i(\pmb x)= \pmb 0$, then $\nabla\Pi_i(\pmb x)=\Pi_{T_{\mathcal{S}_i}(\pmb x)}=\pmb I$.
For this case, the projection $\Pi_i(\cdot)$ is an identity mapping, which implies that the point is already
in $\mathcal{S}_i$. As a result, there is no need to consider this trivial case and only
the nontrivial case of $\nabla f_i(\pmb x)\ne\pmb 0$ is discussed. Denoting the normalized gradients
at $\pmb x^* \in \cap_{i=1}^m \mathcal{S}_i$
\begin{equation}\label{ui:exp}
     \pmb u_i = \frac{\nabla f_i(\pmb x^*)}{\|\nabla f_i(\pmb x^*)\|_2}, ~i=1,\cdots,m
\end{equation}
and according to \eqref{JPi:Tan} and \eqref{Proj:Tan:Si}, the Jacobian matrix of the projection onto the
nonlinear set $\mathcal{S}_i$ is written as
\begin{equation}\label{JPi:exp}
\nabla\Pi_i(\pmb x^*) = \pmb I - \pmb u_i\pmb u_i^\top.
\end{equation}
The Jacobian matrix $\nabla\Pi_i(\pmb x^*)$ is the projection matrix onto
the kernel space of $\pmb u_i^\top$, which is the orthogonal
complement of the range space of $\pmb u_i$ and is denoted as $\mathrm{ker}(\pmb u_i^\top)$.
In addition, the non-expansiveness of $\nabla\Pi_i(\pmb x^*)$ implies
\begin{equation}\label{nonexp:JM:Pi}
     \|\nabla\Pi_i(\pmb x^*)\pmb v\|_2 \le \|\pmb v\|_2
\end{equation}
for any vector $\pmb v$ since $\nabla\Pi_i(\pmb x^*)$ is a projection
onto a linear subspace, which is a convex set. Setting $\pmb z_1=\pmb v$
and $\pmb z_2=\pmb 0$ in \eqref{Nonexpan:Convex} and exploiting
$\nabla\Pi_i(\pmb x^*)\pmb 0=\pmb 0$, we get \eqref{nonexp:JM:Pi}.

The Jacobian matrix of the mapping $\Phi(\cdot)$ at $\pmb x^*$ is
\begin{equation}\label{JM:Phi:Def}
    \nabla\Phi(\pmb x^*) = \left.\frac{\partial \Phi(\pmb z)}
    {\partial\pmb z^\top}\right|_{\pmb z=\pmb x^*}
    \in \mathbb{R}^{n\times n}.
\end{equation}
The spectral norm of $\nabla\Phi(\pmb x^*)$, which is defined as
\begin{equation}\label{}
    \|\nabla\Phi(\pmb x^*)\|_2
    = \sup_{\pmb v\ne\pmb 0} \frac{\|\nabla\Phi(\pmb x^*)\pmb v\|_2}{\|\pmb v\|_2}
\end{equation}
is crucial to the convergence analysis, as stated in the following theorem.
\begin{thm}\label{SpecNorm:Less1}
The spectral norm of the Jacobian matrix of the mappings $\Phi_\mathrm{MP}(\cdot)$
$\Phi_\mathrm{CP}(\cdot)$ and $\Phi_\mathrm{RPP}(\cdot)$ at $\pmb x^*$ is less
than $1$, i.e.,
\begin{equation}\label{}
    \|\nabla\Phi_\mathrm{MP}(\pmb x^*)\|_2 < 1,~
    \|\nabla\Phi_\mathrm{CP}(\pmb x^*)\|_2 < 1~\mathrm{and}~
    \|\nabla\Phi_\mathrm{RPP}(\pmb x^*)\|_2 < 1
\end{equation}
if there are $n$ columns of the matrix
\begin{equation}\label{U:ui}
     \pmb U = [\pmb u_1,\cdots,\pmb u_m]\in\mathbb{R}^{n\times m}
\end{equation}
are linearly independent, or equivalently, $\mathrm{rank}(\pmb U)=n$, where
$\pmb u_i = \nabla f_i(\pmb x^*)/\|\nabla f_i(\pmb x^*)\|_2$.
\end{thm}

\emph{Proof.} We first prove the case of MP. By \eqref{JPi:Tan} of
Lemma \ref{Jacob:Proj} and \eqref{JPi:exp} , it follows that
\begin{equation}\label{}
    \nabla\Phi_\mathrm{MP}(\pmb x^*) = \frac{1}{m}\sum_{i=1}^m\nabla\Pi_i(\pmb x^*)
    =\pmb I - \frac{1}{m}\sum_{i=1}^m \pmb u_i\pmb u_i^\top.
\end{equation}
Since $\nabla\Pi_i(\pmb x^*)$ is symmetric and idempotent, $\nabla\Phi_\mathrm{MP}(\pmb x^*)$ is also
symmetric and idempotent. That is, $\nabla\Phi_\mathrm{MP}(\pmb x^*)^\top=\nabla\Phi_\mathrm{MP}(\pmb x^*)$
and $\nabla\Phi_\mathrm{MP}(\pmb x^*)^2=\nabla\Phi_\mathrm{MP}(\pmb x^*)$, which results in
\begin{equation}\label{JM:Phi:Comp}
\begin{aligned}
    \|\nabla\Phi_\mathrm{MP}(\pmb x^*)\pmb v\|_2^2
    &=\pmb v^\top\nabla\Phi_\mathrm{MP}(\pmb x^*)^\top\nabla\Phi_\mathrm{MP}(\pmb x^*)\pmb v\\
    &=\pmb v^\top\nabla\Phi_\mathrm{MP}(\pmb x^*)\pmb v\\
    &=\|\pmb v\|_2^2 - \frac{1}{m}\sum_{i=1}^m \left(\pmb u_i^\top\pmb v\right)^2.
\end{aligned}
\end{equation}
Since $\pmb U$ is of full column rank, i.e., $\mathrm{rank}(\pmb U)=n$,
we have $\mathrm{ker}(\pmb U^\top)=\{\pmb 0\}$ due to
$\mathrm{rank}(\pmb U)+\mathrm{dim}(\mathrm{ker}(\pmb U^\top))=n$, where
$\mathrm{dim}(\cdot)$ denotes the dimension. Therefore, the linear system
of equations
\begin{equation}\label{}
    \pmb U^\top\pmb v =
    \left[ {\begin{array}{*{20}{c}}
   {\pmb u_1^\top\pmb v}  \\
   {\vdots}  \\
   {\pmb u_m^\top\pmb v}  \\
\end{array}} \right]= \pmb 0
\end{equation}
only has zero solution. As a result, it is impossible for
$\pmb u_i^\top\pmb v=0$ for all $i=1,\cdots,m$ if $\pmb v\ne\pmb 0$.
In other words, $\pmb u_i^\top\pmb v\ne0$ holds for one or several or
even all $i\in\{1,\cdots,m\}$, which implies
$\|\nabla\Phi_\mathrm{MP}(\pmb x^*)\pmb v\|_2^2<\|\pmb v\|_2^2$ for $\pmb v\ne\pmb 0$ and hence
\begin{equation}\label{}
    \|\nabla\Phi_\mathrm{MP}(\pmb x^*)\|_2
    =\sup_{\pmb v\ne \pmb 0}\frac{\|\nabla\Phi_\mathrm{MP}(\pmb x^*)\pmb v\|_2}
    {\|\pmb v\|_2}<1.
\end{equation}

For the case of CP, the mapping $\Phi_\mathrm{CP}(\cdot)$ of \eqref{Phi:CP}
is a composite function that has $m$ layers. We write
$\Phi_\mathrm{CP}(\pmb z) = \Pi_m\Pi_{m-1}\cdots\Pi_2\Pi_1(\pmb z)$ layer-by-layer:
\begin{equation}\label{Recursion:vi}
   \pmb v_1 = \Pi_1(\pmb z), \pmb v_2 = \Pi_2(\pmb v_1), \cdots, \pmb v_m = \Pi_m(\pmb v_{m-1}),
   \Phi_\mathrm{CP}(\pmb z) =\pmb v_m
\end{equation}
which can compactly be expressed by the recursion
$\pmb v_i = \Pi_i(\pmb v_{i-1})$ with $\pmb v_0=\pmb z$.
By the multivariate chain rule, the Jacobian of
$\Phi_\mathrm{CP}(\cdot)$ is computed as
\begin{equation}\label{JM:Phi:CP}
\begin{aligned}
    \nabla\Phi_\mathrm{CP}(\pmb z) &= \frac{\partial \Phi_\mathrm{CP}(\pmb z)}{\partial\pmb z^\top}\\
    &=\frac{\partial \pmb v_m}{\partial\pmb v_{m-1}^\top}
    \frac{\partial \pmb v_{m-1}}{\partial\pmb v_{m-2}^\top}
    \cdots
    \frac{\partial \pmb v_2}{\partial\pmb v_1^\top}
    \frac{\partial \pmb v_1}{\partial\pmb z^\top}\\
    &=\nabla\Pi_m(\pmb v_{m-1})\nabla\Pi_{m-1}(\pmb v_{m-2})
    \cdots\nabla\Pi_2(\pmb v_1)\nabla\Pi_1(\pmb z).
\end{aligned}
\end{equation}
Note that we have $\Pi_i(\pmb x^*)=\pmb x^*$ for all $i\in\{1,\cdots,m\}$
due to $\pmb x^*\in \cap_{i=1}^m \mathcal{S}_i$. Based on the recursion
$\pmb v_i = \Pi_i(\pmb v_{i-1})$ in \eqref{Recursion:vi}, we conclude that $\pmb v_i =\pmb x^*$
for $i=1,\cdots,m$ if $\pmb z=\pmb x^*$. It follows from \eqref{JM:Phi:Def}
and \eqref{JM:Phi:CP} that
\begin{equation}\label{}
    \nabla\Phi_\mathrm{CP}(\pmb x^*) =
    \prod\limits_{i=1}^{m} \nabla\Pi_i(\pmb x^*).
\end{equation}
Because all $\{\nabla\Pi_i(\pmb x^*)\}_{i=1}^m$ are non-expansive, i.e.,
$\|\nabla\Pi_i(\pmb x^*)\pmb v\|_2\le \|\pmb v\|_2$ holds for $i=1,\cdots,m$,
$\nabla\Phi_\mathrm{CP}(\pmb x^*)$ is also non-expansive, implying
$\|\nabla\Phi_\mathrm{CP}(\pmb x^*)\pmb v\|_2\le \|\pmb v\|_2$.
Now we prove that it is impossible that $\|\nabla\Phi_\mathrm{CP}(\pmb x^*)\pmb v\|_2=\|\pmb v\|_2$
for any $\pmb v\ne \pmb 0$ by contradiction.
Suppose $\|\nabla\Phi_\mathrm{CP}(\pmb x^*)\pmb v\|_2=\|\pmb v\|_2$, it requires
\begin{equation}\label{Piv:v}
   \|\nabla\Pi_i(\pmb x^*)\pmb v\|_2 = \|\pmb v\|_2,~\mathrm{for}~ i=1,\cdots,m.
\end{equation}
The reason is explained as follows. The norm of a vector will be reduced or kept the same after the operator
$\nabla\Pi_i(\pmb x^*)$ due to the non-expansiveness. If the norm keeps the same after $m$
operations, it requires that the norm keeps the same at every operation. Noting that
$\nabla\Pi_i(\pmb x^*)$ is symmetric and idempotent and according to \eqref{Piv:v},
we have
\begin{equation}\label{Piv:Eq:v}
\begin{aligned}
   \|\nabla\Pi_i(\pmb x^*)\pmb v\|_2^2
    &= \pmb v^\top \nabla\Pi_i(\pmb x^*)^\top\nabla\Pi_i(\pmb x^*)\pmb v \\
   &=\pmb v^\top\nabla\Pi_i(\pmb x^*)\pmb v \\
   &=\pmb v^\top(\pmb I - \pmb u_i\pmb u_i^\top)\pmb v \\
   &= \|\pmb v\|_2^2 - (\pmb u_i^\top\pmb v)^2
\end{aligned}
\end{equation}
for $i=1,\cdots,m$. Comparing \eqref{Piv:v} and \eqref{Piv:Eq:v}
yields $\{\pmb u_i^\top\pmb v=0\}_{i=1}^m$, i.e.,
$\pmb U^\top\pmb v=\pmb 0$. This indicates that the nonzero vector
$\pmb v\in\mathrm{ker}(\pmb U^\top)$ and $\mathrm{dim}(\mathrm{ker}(\pmb U^\top))\ge 1$.
It follows that
\begin{equation}\label{}
    \mathrm{rank}(\pmb U)=n-\mathrm{dim}(\mathrm{ker}(\pmb U^\top))\le n-1
\end{equation}
which contradicts with $\mathrm{rank}(\pmb U)=n$.
Therefore, the case of $\|\nabla\Phi_\mathrm{CP}(\pmb x^*)\pmb v\|_2=\|\pmb v\|_2$ is impossible and
we have $\|\nabla\Phi_\mathrm{CP}(\pmb x^*)\pmb v\|_2<\|\pmb v\|_2$ for any $\pmb v\ne \pmb 0$,
which implies
\begin{equation}\label{}
    \|\nabla\Phi_\mathrm{CP}(\pmb x^*)\|_2
    =\sup_{\pmb v\ne \pmb 0}\frac{\|\nabla\Phi_\mathrm{CP}(\pmb x^*)\pmb v\|_2}
    {\|\pmb v\|_2}<1.
\end{equation}
The proof for the case of the RPP is similar to that of the CP because
\begin{equation}\label{JM:Phi:RPP}
    \nabla\Phi_\mathrm{RPP}(\pmb x^*)
    =\nabla\Pi_{i_m}(\pmb x^*)\nabla\Pi_{i_{m-1}}(\pmb x^*)
    \cdots\nabla\Pi_{i_2}(\pmb x^*)\nabla\Pi_{i_1}(\pmb x^*)
\end{equation}
has a similar expression to $\nabla\Phi_\mathrm{CP}(\pmb x^*)$.
Following the same steps in the proof of the CP, it is easy to obtain
\begin{equation}\label{}
    \|\nabla\Phi_\mathrm{RPP}(\pmb x^*)\|_2<1
\end{equation}
which completes the proof. \hfill $\square$

\subsection{Proof of Theorem \ref{Thm:Converge:CP}}

We first present the theoretical proof of the convergence of the MP, CP and RPP.
The convergence in expectation of the RP and NRP will be provided in the next subsection
since different techniques are adopted for the case of RP and NRP. The proof of
Theorem \ref{Thm:Converge:CP} is as follows.

\emph{Proof of Theorem \ref{Thm:Converge:CP}.}
We begin with the fixed-point iteration of $\pmb x^{k+1} = \Phi(\pmb x^k),~k=0,1,\cdots$, where $\Phi(\cdot)$
can take forms of $\Phi_{\mathrm{MP}}(\cdot)$, $\Phi_{\mathrm{CP}}(\cdot)$ and $\Phi_{\mathrm{RP}}(\cdot)$.
The mapping $\Phi(\cdot)$ is differentiable at $\pmb x^*$ since we have assumed
that all the projections are differentiable at $\pmb x^*$. The differentiability
of $\Phi(\cdot)$ at $\pmb x^*$ allows the following first-order Taylor series
expansion at $\pmb x^*$
\begin{equation}\label{Taylor:Phi}
    \Phi(\pmb x)= \Phi(\pmb x^*) + \nabla \Phi(\pmb x^*)(\pmb x-\pmb x^*)
    + o(\|\pmb x-\pmb x^*\|_2)
\end{equation}
where $ \nabla\Phi(\pmb x^*)$ is the Jacobian matrix of $\Phi(\cdot)$ at $\pmb x^*$,
which is defined in \eqref{JM:Phi:Def}, and $o(\cdot)$ denotes the
higher-order infinitesimal. It follows from \eqref{Taylor:Phi} that
\begin{equation}\label{}
    \lim_{\pmb x \rightarrow \pmb x^*}
    \frac{\|\Phi(\pmb x)-\Phi(\pmb x^*)
    -\nabla \Phi(\pmb x^*)(\pmb x-\pmb x^*)\|_2}
    {\|\pmb x-\pmb x^*\|_2}=0.
\end{equation}
By the $(\epsilon,\delta)$-definition of the limit of a function \cite{Royden}, there
exists a neighborhood around $\pmb x^*$ with radius of $\delta$, which is
denoted as $\mathcal{B}_\delta(\pmb x^*) = \left\{\pmb x\left|\|\pmb x
- \pmb x^*\| < \delta\right\}\right.$, such that
$\forall \pmb x\in \mathcal{B}_\delta(\pmb x^*)$, it has
\begin{equation}\label{deriv:lim}
    \frac{\|\Phi(\pmb x)-\Phi(\pmb x^*)
    -\nabla \Phi(\pmb x^*)(\pmb x-\pmb x^*)\|_2}
    {\|\pmb x-\pmb x^*\|_2}<\epsilon
\end{equation}
where $\epsilon>0$ can be arbitrarily small. The radius of the neighborhood
$\delta$ depends on $\epsilon$. Based on the fact that $\pmb x^*$ is a
fixed-point of $\Phi(\cdot)$, i.e., $\pmb x^*=\Phi(\pmb x^*)$, we compute
\begin{align}\nonumber
\|\Phi(\pmb x)-\pmb x^*\|_2
&= \|\Phi(\pmb x)-\Phi(\pmb x^*)\|_2\\ \nonumber
&= \|\Phi(\pmb x)-\Phi(\pmb x^*)-\nabla \Phi(\pmb x^*)(\pmb x-\pmb x^*)
+\nabla \Phi(\pmb x^*)(\pmb x-\pmb x^*)\|_2\\ \label{Conv:Ineq1}
&\le \|\Phi(\pmb x)-\Phi(\pmb x^*)-\nabla \Phi(\pmb x^*)(\pmb x-\pmb x^*)\|_2
+ \|\nabla \Phi(\pmb x^*)(\pmb x-\pmb x^*)\|_2 \\ \label{Conv:Ineq2}
&<\epsilon \|\pmb x-\pmb x^*\|_2 + \|\nabla \Phi(\pmb x^*)\|_2\|\pmb x-\pmb x^*\|_2\\ \nonumber
&= (\epsilon+\|\nabla \Phi(\pmb x^*)\|_2)\|\pmb x-\pmb x^*\|_2\\ \nonumber
&= \gamma\|\pmb x-\pmb x^*\|_2
\end{align}
where \eqref{Conv:Ineq1} follows from the triangle inequality and \eqref{Conv:Ineq2} is
due to \eqref{deriv:lim} and the spectral norm inequality
$\|\nabla \Phi(\pmb x^*)(\pmb x-\pmb x^*)\|_2\le\|\nabla \Phi(\pmb x^*)\|_2\|\pmb x-\pmb x^*\|_2$.
Because $\epsilon>0$ can be arbitrarily small and $\|\nabla \Phi(\pmb x^*)\|_2<1$
holds true when $\Phi(\cdot)$ takes forms of $\Phi_\mathrm{MP}(\cdot)$,
$\Phi_\mathrm{CP}(\cdot)$ and $\Phi_\mathrm{RPP}(\cdot)$ by Theorem
\ref{SpecNorm:Less1}, we see that
\begin{equation}\label{}
    \gamma = \epsilon + \|\nabla \Phi(\pmb x^*)\|_2<1
\end{equation}
is guaranteed. If $\pmb x^k\in\mathcal{B}_\delta(\pmb x^*)$, it follows that
\begin{equation}
\|\pmb x^{k+1}-\pmb x^*\|_2
=\|\Phi(\pmb x^k)-\pmb x^*\|_2 <\gamma\|\pmb x^k-\pmb x^*\|_2
\end{equation}
which implies that the iteration is contracted and
$\pmb x^{k+1}\in\mathcal{B}_\delta(\pmb x^*)$. By induction,
the convergence with a linear rate
\begin{equation}
\|\pmb x^k - \pmb x^*\|_2 < \gamma^k\|\pmb x^0-\pmb x^*\|_2
\end{equation}
is obtained, which completes the proof. \hfill $\square$

\emph{Proof of Corollary \ref{Thm:Rate:MP}.} According to \eqref{JM:Phi:Comp} and rewriting
$\mathop{\sum}\limits_{i=1}^m \left(\pmb u_i^\top\pmb v\right)^2
=\|\pmb U^\top \pmb v\|_2^2$,
we obtain
\begin{equation}\label{}
    \|\nabla\Phi_\mathrm{MP}(\pmb x^*)\|_2^2
    = \max_{\|\pmb v\|_2=1}\|\nabla\Phi_\mathrm{MP}(\pmb x^*) \pmb v\|_2^2
    =1 - \frac{1}{m}\min_{\|\pmb v\|_2=1}\|\pmb U^\top \pmb v\|_2^2.
\end{equation}
By Courant-Fischer theorem \cite{Golub}, $\mathop{\min}\limits_{\|\pmb v\|_2=1}
\|\pmb U^\top \pmb v\|_2^2=\lambda_{\min}(\pmb U\pmb U^\top)
=\sigma_{\min}^2(\pmb U)$ with $\lambda_{\min}(\cdot)$ being
the minimum eigenvalue of a matrix. \hfill $\square$

\subsection{Proof of Theorem \ref{Thm:Rate:RP}}

\emph{Proof of Theorem \ref{Thm:Rate:RP}.} The iteration of RP and NRP is $\pmb x^{k+1} = \Pi_{i_k}(\pmb x^k)$ with
$i_k$ being uniformly or non-uniformly sampled from $\{1,\cdots,m\}$. Since $\Pi_{i_k}(\cdot)$ is
differentiable at $\pmb x^*$, taking the same steps as in
\eqref{Taylor:Phi}--\eqref{deriv:lim}, it concludes that
\begin{equation}\label{RP:Deriv:Ineq}
    \|\Pi_{i_k}(\pmb x)-\Pi_{i_k}(\pmb x^*)
    -\nabla \Pi_{i_k}(\pmb x^*)(\pmb x-\pmb x^*)\|_2<\epsilon\|\pmb x-\pmb x^*\|_2
\end{equation}
holds true $\forall \pmb x\in \mathcal{B}_\delta(\pmb x^*)$ where
$\epsilon>0$ can be arbitrarily small. Starting with $\pmb x^*=\Pi_{i_k}(\pmb x^*)$,
we derive
\begin{align}\nonumber
\|\Pi_{i_k}(\pmb x)-\pmb x^*\|_2
&= \|\Pi_{i_k}(\pmb x)-\Pi_{i_k}(\pmb x^*)\|_2\\ \label{RP:TriAngleIn}
&\le \|\Pi_{i_k}(\pmb x)-\Pi_{i_k}(\pmb x^*)-\nabla \Pi_{i_k}(\pmb x^*)(\pmb x-\pmb x^*)\|_2
+ \|\nabla \Pi_{i_k}(\pmb x^*)(\pmb x-\pmb x^*)\|_2 \\ \label{RP:Ineq2}
&<\epsilon \|\pmb x-\pmb x^*\|_2 + \|\nabla \Pi_{i_k}(\pmb x^*)(\pmb x-\pmb x^*)\|_2
\end{align}
where \eqref{RP:TriAngleIn} follows from the triangle inequality and \eqref{RP:Ineq2} is
due to \eqref{RP:Deriv:Ineq}. Taking expectation on
the left-hand and right-hand sides of the above inequality yields
\begin{align}\nonumber
\mathbb{E}_{i_k}\left[\|\pmb x^{k+1}-\pmb x^*\|_2\right]
&=\mathbb{E}_{i_k}\left[\|\Pi_{i_k}(\pmb x^k)-\pmb x^*\|_2\right]\\ \nonumber
&<\epsilon \|\pmb x^k-\pmb x^*\|_2 + \mathbb{E}_{i_k}\left[
\|\nabla \Pi_{i_k}(\pmb x^*)(\pmb x^k-\pmb x^*)\|_2\right]\\ \nonumber
&=\epsilon \|\pmb x^k-\pmb x^*\|_2 + \sqrt{\mathbb{E}_{i_k}^2\left[
\|\nabla \Pi_{i_k}(\pmb x^*)(\pmb x^k-\pmb x^*)\|_2\right]}\\ \label{E2:Ineq}
&\le\epsilon \|\pmb x^k-\pmb x^*\|_2 + \sqrt{\mathbb{E}_{i_k}
\left[\|\nabla \Pi_{i_k}(\pmb x^*)(\pmb x^k-\pmb x^*)\|_2^2\right]}
\end{align}
where \eqref{E2:Ineq} is based on the inequality $\mathbb{E}^2[\xi]\le \mathbb{E}[\xi^2]$
with $\xi\in\mathbb{R}$ being a random variable, which can be derived from
$\mathbb{E}[(\xi-\mathbb{E}[\xi])^2]\ge 0$. Exploiting \eqref{JPi:exp},
employing again that $\nabla\Pi_i(\pmb x^*)$ is symmetric and idempotent
and following similar steps in \eqref{Piv:Eq:v}, we obtain
\begin{equation}\label{Pi:x:xstar}
   \|\nabla\Pi_i(\pmb x^*)(\pmb x^k-\pmb x^*)\|_2^2
   =\|\pmb x^k-\pmb x^*\|_2^2 - \left(\pmb u_i^\top(\pmb x^k-\pmb x^*)\right)^2.
\end{equation}
For RP, the uniform distribution over $\{1,\cdots,m\}$ of $i_k$ implies that
\begin{equation}\label{}
\begin{aligned}
   \mathbb{E}_{i_k}\left[\|\nabla\Pi_{i_k}(\pmb x^*)(\pmb x^k-\pmb x^*)\|_2^2\right]
   &=\frac{1}{m}\sum_{i=1}^m \|\nabla \Pi_i(\pmb x^*)(\pmb x^k-\pmb x^*)\|_2^2\\
   &=\|\pmb x^k-\pmb x^*\|_2^2 - \frac{1}{m}\sum_{i=1}^m \left(\pmb u_i^\top(\pmb x^k-\pmb x^*)\right)^2\\
   &=\|\pmb x^k-\pmb x^*\|_2^2 - \frac{1}{m}\|\pmb U^\top(\pmb x^k-\pmb x^*)\|_2^2\\
   &\le\|\pmb x^k-\pmb x^*\|_2^2 - \frac{\sigma_{\min}^2(\pmb U)}{m}\|\pmb x^k-\pmb x^*\|_2^2\\
   &=\left(1 - \frac{\sigma_{\min}^2(\pmb U)}{m}\right)\|\pmb x^k-\pmb x^*\|_2^2.
\end{aligned}
\end{equation}
Plugging the above inequality into \eqref{E2:Ineq} leads to
\begin{equation}
\mathbb{E}_{i_k}\left[\|\pmb x^{k+1}-\pmb x^*\|_2\right]
<\left(\epsilon + \sqrt{1 - \frac{\sigma_{\min}^2(\pmb U)}{m}}\right)\|\pmb x^k-\pmb x^*\|_2
\end{equation}
where the convergence rate
\begin{equation}
\gamma_{\mathrm{RP}} = \epsilon + \sqrt{1 - \frac{\sigma_{\min}^2(\pmb U)}{m}}<1
\end{equation}
since $\sigma_{\min}^2(\pmb U)>0$ due to $\mathrm{rank}(\pmb U)=n$
and $\epsilon>0$ can be arbitrarily small. Thus, the iteration of RP
is contracted and convergent in expectation. By induction, we obtain the
following convergence with a linear rate
\begin{equation}\label{}
    \mathbb{E}_{i_k}\left[\|\pmb x^{k+1} - \pmb x^*\|_2\right]
    < \gamma_\mathrm{RP}^{k+1}\|\pmb x^0 - \pmb x^*\|_2.
\end{equation}
The asymptotic convergence rate of RP is
\begin{equation}\label{}
    \mathop{\lim}\limits_{k\rightarrow\infty}\gamma_\mathrm{RP}
    =\sqrt{1 - \frac{1}{\kappa^2(\pmb U)}}
\end{equation}
where $\kappa(\pmb U)=\frac{\|\pmb U\|_\mathrm{F}}{\sigma_{\min}(\pmb U)}$
is the condition number of $\pmb U$.

Observing that the nonuniform distribution of \eqref{NRP:rule} for NRP
can be rewritten as
\begin{equation}\label{NRP:rule2}
  \frac{\|\nabla f_i(\pmb x^*)\|_2^2}{\|\pmb G\|_\mathrm{F}^2},
  \quad i=1,\cdots,m
\end{equation}
we have for NRP that
\begin{equation}\label{}
\begin{aligned}
   \mathbb{E}_{i_k}\left[\|\nabla\Pi_{i_k}(\pmb x^*)(\pmb x^k-\pmb x^*)\|_2^2\right]
   &=\sum_{i=1}^m \frac{\|\nabla f_i(\pmb x^*)\|_2^2}{\|\pmb G\|_\mathrm{F}^2}
   \|\nabla \Pi_i(\pmb x^*)(\pmb x^k-\pmb x^*)\|_2^2\\
   &=\|\pmb x^k-\pmb x^*\|_2^2
   -\frac{1}{\|\pmb G\|_\mathrm{F}^2}\sum_{i=1}^m\left(\nabla f_i(\pmb x^*)^\top(\pmb x^k-\pmb x^*)\right)^2\\
   &=\|\pmb x^k-\pmb x^*\|_2^2 - \frac{1}{\|\pmb G\|_\mathrm{F}^2}\|\pmb G^\top(\pmb x^k-\pmb x^*)\|_2^2\\
   &\le\|\pmb x^k-\pmb x^*\|_2^2 - \frac{\sigma_{\min}^2(\pmb G)}{\|\pmb G\|_\mathrm{F}^2}\|\pmb x^k-\pmb x^*\|_2^2\\
   &=\left(1 - \frac{\sigma_{\min}^2(\pmb G)}{\|\pmb G\|_\mathrm{F}^2}\right)\|\pmb x^k-\pmb x^*\|_2^2\\
   &=\left(1 - \frac{1}{\kappa^2(\pmb G)}\right)\|\pmb x^k-\pmb x^*\|_2^2
\end{aligned}
\end{equation}
where we have used $\mathop{\sum}\limits_{i=1}^m \frac{\|\nabla f_i(\pmb x^*)\|_2^2}{\|\pmb G\|_\mathrm{F}^2}=1$
and $\|\nabla f_i(\pmb x^*)\|_2\pmb u_i=\nabla f_i(\pmb x^*)$.
The remaining steps are the same as those of the RP, which completes the proof. \hfill $\square$

\begin{remark}
Comparing \eqref{AsymRate:MP} and \eqref{AsymRate:RP}, we see that the MP and RP have
the same asymptotic convergence rate. However, the MP requires $m$ projections in one
iteration while the RP just needs one. Therefore, they are different and the RP is $m$
times faster than the MP.
\end{remark}

\subsection{Proof of Theorem \ref{Thm:Rate:GP}}

We need the following lemma to analyze the convergence of the GP and NGP.

\begin{lem}[Gradient Representation of Greedy Rules]\label{Lemma:Greedy:Gradient}
There exists a neighborhood centered at $\pmb x^*$ with radius $\delta'$
\begin{equation}\label{}
    \mathcal{B}_{\delta'}(\pmb x^*) =
    \left\{\pmb x\left|\|\pmb x - \pmb x^*\|_2 < \delta'\right\} \right.
\end{equation}
such that if $\pmb x^k\in \mathcal{B}_{\delta'}(\pmb x^*)$, then
the GP rule of \eqref{GP:rule} is equivalent to
\begin{equation}\label{GP:rule2}
    i_k= \arg\max_{1\le i\le m}\left|\nabla f_i(\pmb x^*)^\top(\pmb x^k-\pmb x^*)\right|
\end{equation}
and the NGP rule of \eqref{NGP:rule} amounts to
\begin{equation}\label{NGP:rule2}
    i_k= \arg\max_{1\le i\le m}\left|\pmb u_i^\top(\pmb x^k-\pmb x^*)\right|.
\end{equation}
\end{lem}
\emph{Proof.} Using $f_i(\pmb x^*)=0$, we obtain the first-order Taylor
series expansion of $f_i(\pmb x)$ at $\pmb x^*$
\begin{equation}\label{Taylor:fi}
    f_i(\pmb x)= \nabla f_i(\pmb x^*)^\top(\pmb x-\pmb x^*)
    + o(\|\pmb x-\pmb x^*\|_2)
\end{equation}
which indicates that $\nabla f_i(\pmb x^*)^\top(\pmb x-\pmb x^*)$ is a
good approximation of $f_i(\pmb x)$ if $\pmb x$ is close enough to $\pmb x^*$.
There exists $\delta'>0$ when $\pmb x^k\in \mathcal{B}_{\delta'}(\pmb x^*)$
such that the approximation $f_i(\pmb x^k)\approx \nabla f_i(\pmb x^*)^\top(\pmb x^k-\pmb x^*)$
is accurate enough for all $i=1,\cdots,m$. Therefore, if $\delta'>0$ is small enough, then
the accurate first-order Taylor approximation guarantees that
$\arg\mathop{\max}\limits_{1\le i\le m}|f_i(\pmb x^k)|$ is equivalent to
$\arg\mathop{\max}\limits_{1\le i\le m}\left|\nabla f_i(\pmb x^*)^\top(\pmb x^k-\pmb x^*)\right|$
and $\arg\mathop{\max}\limits_{1\le i\le m}\frac{\left|f_i(\pmb x^k)\right|}{\|\nabla f_i(\pmb x^*)\|_2}$
amounts to $\arg\mathop{\max}\limits_{1\le i\le m}\left|\pmb u_i^\top(\pmb x^k-\pmb x^*)\right|$,
recalling that $\pmb u_i=\frac{\nabla f_i(\pmb x^*)}{\|\nabla f_i(\pmb x^*)\|_2}$.
The proof is complete. \hfill $\square$

Now we are ready to present the proof of Theorem \ref{Thm:Rate:GP}.

\emph{Proof of Theorem \ref{Thm:Rate:GP}.} It follows from \eqref{RP:Ineq2} and \eqref{Pi:x:xstar} that
\begin{equation}\label{GP:Basic:Ineq}
\begin{aligned}
\|\pmb x^{k+1}-\pmb x^*\|_2
&=\|\Pi_{i_k}(\pmb x^k)-\pmb x^*\|_2\\
&<\epsilon \|\pmb x^k-\pmb x^*\|_2 + \|\nabla \Pi_{i_k}(\pmb x^*)(\pmb x^k-\pmb x^*)\|_2\\
&=\epsilon \|\pmb x^k-\pmb x^*\|_2
+ \sqrt{\|\pmb x^k-\pmb x^*\|_2^2
- \left(\pmb u_{i_k}^\top(\pmb x^k-\pmb x^*)\right)^2}.
\end{aligned}
\end{equation}
By \eqref{NGP:rule2} in Lemma \ref{Lemma:Greedy:Gradient}, we have for the NGP that
\begin{equation}\label{ui:hoff:inf}
\begin{aligned}
    \left(\pmb u_{i_k}^\top(\pmb x^k-\pmb x^*)\right)^2
    &=\max_{1\le i\le m}\left(\pmb u_i^\top(\pmb x^k-\pmb x^*)\right)^2\\
    &=\|\pmb U^\top(\pmb x^k-\pmb x^*)\|_\infty^2\\
    &\ge h_\infty^2(\pmb U^\top)\|\pmb x^k-\pmb x^*\|_2^2
\end{aligned}
\end{equation}
where $h_\infty(\pmb U^\top)$ is the Hoffman type
constant \cite{Hoffman} of $\pmb U^\top$, which is defined as
\begin{equation}\label{}
  h_p(\pmb U^\top) = \inf_{\pmb v\ne \pmb 0}\frac{\|\pmb U^\top \pmb v\|_p}{\|\pmb v\|_2}
\end{equation}
with $\|\cdot\|_p$ denoting the $\ell_p$-norm. For $p=2$,
$h_2(\pmb U^\top)=\sigma_{\min}(\pmb U^\top)$ is
the minimum singular value of $\pmb U^\top$. Using
the inequality $\|\pmb a\|_2\le\sqrt{m}\|\pmb a\|_\infty$
with $\pmb a\in\mathbb{R}^m$, it is easy to verify that
\begin{equation}\label{}
  h_\infty(\pmb U^\top)\ge \frac{1}{\sqrt{m}}h_2(\pmb U^\top)
  = \frac{1}{\sqrt{m}}\sigma_{\min}(\pmb U^\top)>0
\end{equation}
if $\mathrm{rank}(\pmb U)=n$. Substituting \eqref{ui:hoff:inf}
into \eqref{GP:Basic:Ineq} yields
\begin{equation}\label{}
\|\pmb x^{k+1}-\pmb x^*\|_2
<\left(\epsilon + \sqrt{1-h_\infty^2(\pmb U^\top)}\right) \|\pmb x^k-\pmb x^*\|_2
\end{equation}
which elicits
\begin{equation}\label{}
\|\pmb x^k-\pmb x^*\|_2
<\gamma_\mathrm{NGP}^k\|\pmb x^0-\pmb x^*\|_2
\end{equation}
where the convergence rate
\begin{equation}\label{}
\gamma_\mathrm{NGP}=\epsilon + \sqrt{1-h_\infty^2(\pmb U^\top)}<1
\end{equation}
since $1-h_\infty^2(\pmb U^\top)<1$ and $\epsilon>0$ can be arbitrarily small.
The asymptotic convergence rate of the NGP is
\begin{equation}\label{}
    \mathop{\lim}\limits_{k\rightarrow\infty}\gamma_\mathrm{NGP}
    =\sqrt{1-h_\infty^2(\pmb U^\top)}.
\end{equation}
It is obvious that the NGP is faster than the RP due to
$h_\infty^2(\pmb U^\top)\ge\frac{1}{m}\sigma_{\min}^2(\pmb U^\top)$.
Note that the worst case where the equality $h_\infty^2(\pmb U^\top)
=\frac{1}{m}\sigma_{\min}^2(\pmb U^\top)$ holds seldom occurs.
Numerical experiment results in Section \ref{Sec:Simulation} show
that the NGP is much faster than the RP in general.
According to \eqref{GP:rule2}, we obtain
\begin{equation}\label{}
    \left(\nabla f_{i_k}(\pmb x^*)^\top(\pmb x^k-\pmb x^*)\right)^2
    =\max_{1\le i\le m}\left(\nabla f_i(\pmb x^*)^\top(\pmb x^k-\pmb x^*)\right)^2
    =\|\pmb G^\top(\pmb x^k-\pmb x^*)\|_\infty^2
\end{equation}
for the GP, which results in
\begin{equation}\label{ui:hoff:GP}
    \begin{aligned}
    \left(\pmb u_{i_k}^\top(\pmb x^k-\pmb x^*)\right)^2
    &=\frac{\left(\nabla f_{i_k}(\pmb x^*)^\top(\pmb x^k-\pmb x^*)\right)^2}{\|\nabla f_{i_k}(\pmb x^*)\|_2^2}\\
    &=\frac{\|\pmb G^\top(\pmb x^k-\pmb x^*)\|_\infty^2}{\|\nabla f_{i_k}(\pmb x^*)\|_2^2}\\
    &\ge \frac{h_\infty^2(\pmb G^\top)}{\|\nabla f_{i_k}(\pmb x^*)\|_2^2}\|\pmb x^k-\pmb x^*\|_2^2.
\end{aligned}
\end{equation}
Plugging \eqref{ui:hoff:GP} into \eqref{GP:Basic:Ineq} leads to
\begin{equation}\label{GP:Relation}
\|\pmb x^{k+1}-\pmb x^*\|_2
<\left(\epsilon + \sqrt{1-\frac{h_\infty^2(\pmb G^\top)}{\|\nabla f_{i_k}(\pmb x^*)\|_2^2}}\right)
\|\pmb x^k-\pmb x^*\|_2.
\end{equation}
Recursively applying \eqref{GP:Relation}, we obtain
\begin{equation}\label{GP:rate}
\|\pmb x^k-\pmb x^*\|_2
<\prod_{j=1}^k\left(\epsilon_j + \sqrt{1-\frac{h_\infty^2(\pmb G^\top)}{\|\nabla f_{i_j}(\pmb x^*)\|_2^2}}\right)
\|\pmb x^0-\pmb x^*\|_2
\end{equation}
where we emphasize that $\epsilon$ depends on the iteration number $j$
and is denoted as $\epsilon_j$. As the iteration progresses,
$\mathop{\lim}\limits_{j\rightarrow\infty}\epsilon_j=0$. This convergence
rate is related to the specific $\|\nabla f_{i_j}(\pmb x^*)\|_2$ with $i_j$
being the selected index at the $j$th iteration. By noticing that the
$\ell_{2,\infty}$-norm of $\pmb G^\top$ is defined as
\begin{equation}\label{}
\left\|\pmb G^\top\right\|_{2,\infty}=\max_{1\le i\le m}\|\nabla f_i(\pmb x^*)\|_2
\end{equation}
we can obtain a looser but more concise bound
\begin{equation}\label{GP:Relation2}
\|\pmb x^{k+1}-\pmb x^*\|_2
<\left(\epsilon + \sqrt{1-\frac{h_\infty^2(\pmb G^\top)}{\left\|\pmb G^\top\right\|_{2,\infty}^2}}\right)
\|\pmb x^k-\pmb x^*\|_2
\end{equation}
which is equivalent to
\begin{equation}\label{GP:rate2}
\|\pmb x^k-\pmb x^*\|_2
<\gamma_\mathrm{GP}^k\|\pmb x^0-\pmb x^*\|_2
\end{equation}
where the convergence rate
\begin{equation}\label{}
\gamma_\mathrm{GP}=\epsilon +
\sqrt{1-\frac{h_\infty^2(\pmb G^\top)}{\left\|\pmb G^\top\right\|_{2,\infty}^2}}<1.
\end{equation}
The asymptotic convergence rate of the GP is
\begin{equation}\label{}
    \mathop{\lim}\limits_{k\rightarrow\infty}\gamma_\mathrm{GP}
    =\sqrt{1-\frac{h_\infty^2(\pmb G^\top)}{\left\|\pmb G^\top\right\|_{2,\infty}^2}}.
\end{equation}
We emphasize that the convergence rate of \eqref{GP:rate} is tighter than
that of \eqref{GP:rate2} while the latter is just more concise. \hfill $\square$

Noting that $\pmb G^\top=\nabla\pmb f(\pmb x^*)$
is the Jacobian matrix of $\pmb f(\cdot)$ at the solution
$\pmb x^*$ and $\pmb U^\top$ is the Jacobian matrix with normalized columns, now we can summarize
from Corollary \ref{Thm:Rate:MP} and Theorems \ref{Thm:Rate:RP} and \ref{Thm:Rate:GP} that the
convergence rates of the variants of the SP depend on the Hoffman constants
of the Jacobian matrix of the nonlinear functions. Due to the linear rate of convergence,
the iteration complexity of all variants of the SP is
$\mathcal{O}(\log(1/\eta))$ to achieve an $\eta$-accuracy solution.

\section{Numerical Results}\label{Sec:Simulation}

The convergence behaviors of the variants of the SP, including the CP, RP, RPP, NRP, GP and NGP are
investigated in solving the phase retrieval and circle equations problems.

\subsection{Results of Phase Retrieval}

In the first simulation, we consider the phase retrieval problem, i.e., solving the system
of elliptic equations $\left|\pmb a_i^\mathrm{H}\pmb x\right|^2 = b_i^2$ with $i=1,\cdots,m$.
See Appendix A for details of this problem. Both $\pmb x\in\mathbb{C}^n$ and $\{\pmb a_i\}$ are randomly generated from a complex standard
i.i.d. Gaussian distribution. We set $n=128$ and $m=5n$. At each iteration, the SP projects
the current solution onto the surface of one ellipsoid according to \eqref{Proj:Ei}.

\begin{figure}
\begin{center}
\includegraphics[width=10cm]{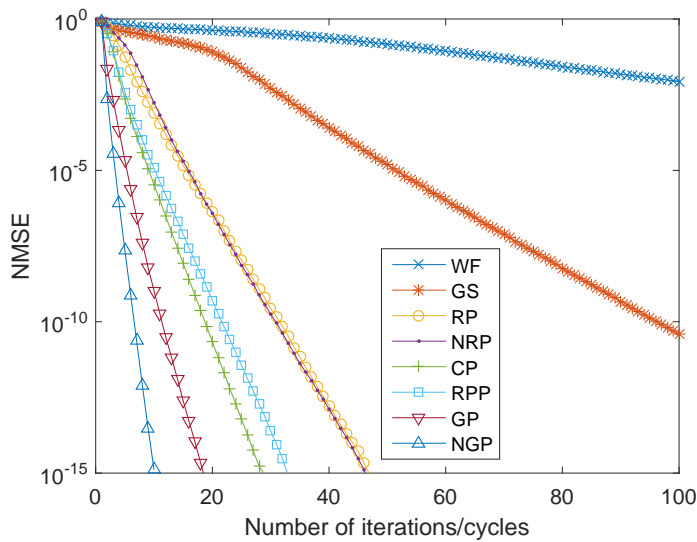}
\caption{NMSE versus number of iterations/cycles.}\label{Fig:PR:iter}
\end{center}
\end{figure}

Two well-known phase retrieval methods, namely, the Wirtinger flow (WF) \cite{WF}
and Gerchberg-Saxton (GS) algorithm \cite{GSA}, are compared with our solvers.
It is fair to compare $m$ iterations (one cycle) for the SP with one WF iteration
because the time complexity of the CP, RP, RPP and NRP per cycle is $\mathcal{O}(mn)$,
which is the same as the WF per iteration. The GP, NGP and GS has a higher
complexity of $\mathcal{O}(m^2n)$ per cycle. But still, we plot the results of
the three methods per cycle for comparison. Since there is an intrinsic phase
ambiguity in phaseless equations, the following normalized mean squared error (NMSE)
with a phase alignment
\begin{equation*}
\mathrm{NMSE}(\pmb x^k) = \frac{\min_{\phi\in[0,2\pi)}
\|\pmb x^k - {\rm e}^{{\rm j}\phi}\pmb x^* \|_2^2}{\|\pmb x^*\|_2^2}
\end{equation*}
is taken as the performance index, where $\pmb x^*$ is the original signal (true solution).
This index reflects the speed of convergence to the original signal.
Figure~\ref{Fig:PR:iter} shows the NMSE versus the number of iterations/cycles. We observe that
all methods converge to the original signal (global solution) $\pmb x^*$ at a linear rate. The six SP solvers
converge much faster than the WF and GS. The convergence speeds of the two greedy solvers
are the fastest among them.

\subsection{Results of Solving Circle Equations}

In the second simulation, we consider solving the circle equations
$\|\pmb x - \pmb c_i\|_2^2 = r_i^2$, $i=1,\cdots,m$. See Appendix A
for its applications to range measurement localization.
In our experiment, the centers $\{\pmb c_i\}\in \mathbb{R}^n$ and the true solution
$\pmb x^*\in \mathbb{R}^n$ are randomly generated from the standard i.i.d.
Gaussian distribution. Once $\{\pmb c_i\}$ and
$\pmb x^*$ are generated, the radiuses $\{r_i\}$ can be computed.
We set $n=100$ and $m=400$. At each iteration, the SP projects
the current solution onto one sphere according to \eqref{Proj:Ci}.

\begin{figure}
\begin{center}
\includegraphics[width=10cm]{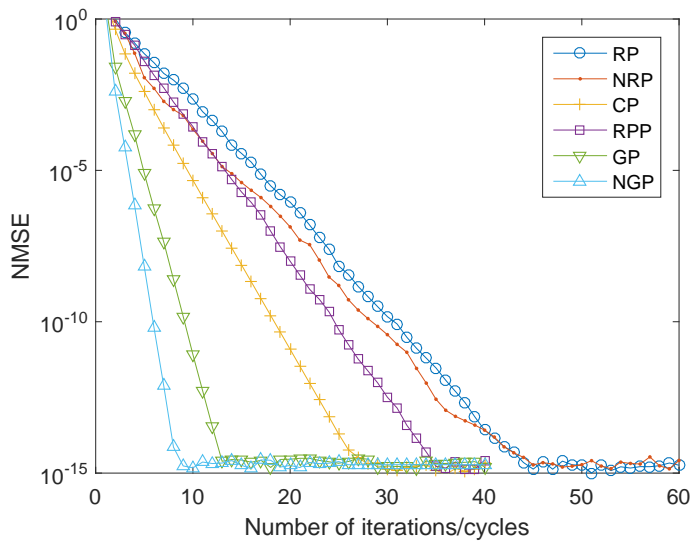}
\caption{NMSE versus number of iterations/cycles.}\label{Fig:CE:iter}
\end{center}
\end{figure}

The following NMSE
\begin{equation*}
\mathrm{NMSE}(\pmb x^k) = \frac{\|\pmb x^k - \pmb x^* \|_2^2}{\|\pmb x^*\|_2^2}
\end{equation*}
is adopted as the performance index, where $\pmb x^*$ is the true solution.
This index reflects the speed of convergence to the true solution.
Figure~\ref{Fig:CE:iter} displays the NMSE versus the number of iterations/cycles
of the six solvers. We observe that all methods converge to the true solution
$\pmb x^*$ at a linear rate. The six SP solvers converge very fast.
Again, the convergence rates of the two greedy solvers are the fastest.

\section*{A\quad Important Nonlinear Equations/Inequalites}

Some nonlinear equations/inequalites frequently encountered in
engineering and science are listed in the following.
\begin{itemize}
  \item Circle equations:
  \begin{equation}\label{Cir:Eq}
    f_i(\pmb x)=\|\pmb x - \pmb c_i\|_2^2 - r_i^2 = 0, \quad i=1,\cdots,m
  \end{equation}
  where $\pmb c_i\in \mathbb{R}^n$ and $r_i>0$ are the center and
  radius of the $i$th ball, respectively. Solving circle equations has
  important applications in range measurement based source localization.
  Source localization aims at determining the position of a source $\pmb x$.
  There are $m$ receiving sensors, whose positions are known as $\{\pmb c_i\}_{i=1}^m$,
  measuring the Euclidean distances $\{r_i\}_{i=1}^m$ from the source to the receivers.
  The solution of the circle equations \eqref{Cir:Eq} gives the source position.
  Since the sphere constraints are nonconvex, solving circle equations refers to
  a nonconvex feasibility problem.
  \item Phase retrieval: phase retrieval refers to recovering a complex-valued
  signal $\pmb x\in\mathbb{C}^n$ from $m$ squared magnitude-only
  measurements
\begin{equation}\label{PR:Model}
      b_i^2 = |\pmb a_i^\mathrm{H}\pmb x|^2, \quad i = 1,\cdots,m
\end{equation}
where $\pmb a_i\in\mathbb{C}^n$ are known sampling vectors and $b_i\in\mathbb{R}_+$.
Finding a solution of \eqref{PR:Model} implies to solve a
system of quadratic equations
\begin{equation}\label{Phase:Retr}
f_i(\pmb x)=\pmb x^\mathrm{H}\pmb a_i\pmb a_i^\mathrm{H}\pmb x - b_i^2 = 0.
\end{equation}
Geometrically, \eqref{Phase:Retr} is the surface of an (degenerate) ellipsoid, which is nonconvex.
\item General quadratic equations:
  \begin{equation}\label{GQ:Eq}
    f_i(\pmb x)=\pmb x^\top\pmb A_i\pmb x - 2\pmb c_i^\top\pmb x  + b_i = 0.
  \end{equation}
  where $\pmb A_i\in\mathbb{R}^{n \times n}$ may not be positive semi-definite,
  $\pmb c_i\in \mathbb{R}^n$ and $b_i\in\mathbb{R}$.
  \item Sparsity constraint: in compressed sensing, the signal is often sparse in a specific basis.
  The following (quasi-) $\ell_p$-norm constraint is used to promote the sparsity
  \begin{equation}\label{LowRank:Cons}
    \|\pmb x\|_p \le s
  \end{equation}
  where $0\le p\le 1$ and $s\in\mathbb{R}_+$ controls the sparsity.
  The $\ell_1$-ball constraint is convex for $p=1$ while it is nonconvex for $0\le p<1$.
\item Rank constraint: in matrix completion and low-rank approximation, the matrices often
have the low-rank property, yielding the following low-rank constraint
  \begin{equation}\label{LowRank:Cons}
    \mathrm{rank}(\pmb X)\le r
  \end{equation}
  where the variable $\pmb X\in\mathbb{R}^{m \times n}$ is in the form of a matrix
  and $r\le\min(m,n)$ is the desired rank. The rank constraint is nonconvex.
\end{itemize}

\section*{B\quad Norms of the Gradients}

In many practical applications, the norms of the gradients $\{\|\nabla f_i(\pmb x^*)\|_2\}$
can be computed without knowing the solution $\pmb x^*$. Some examples include:
\begin{itemize}
  \item Linear equation: $f_i(\pmb x) = \pmb a_i^\top\pmb x - b_i$
  with $\|\nabla f_i(\pmb x^*)\|_2 = \|\pmb a_i\|_2$.
  \item Circle equation: $f_i(\pmb x) =\|\pmb x - \pmb c_i\|_2^2 - r_i^2$
  with $\|\nabla f_i(\pmb x^*)\|_2 = 2r_i$.
  \item Unsquared circle equation: $f_i(\pmb x) =\|\pmb x - \pmb c_i\|_2 - r_i$
  with $\|\nabla f_i(\pmb x^*)\|_2 = 1$. Note that $\nabla f_i(\pmb x)=
  \frac{\pmb x - \pmb c_i}{\|\pmb x - \pmb c_i\|_2}$ and it is
  well-defined at $\pmb x^*$ due to $\pmb x^*\ne\pmb c_i$.
  \item Elliptic equation: $f_i(\pmb x)=\pmb x^\mathrm{H}\pmb a_i\pmb a_i^\mathrm{H}\pmb x - b_i^2$
  with $\|\nabla f_i(\pmb x^*)\|_2 = b_i\|\pmb a_i\|_2$.
  \item Unsquared elliptic equation: $f_i(\pmb x)=|\pmb a_i^\mathrm{H}\pmb x| - b_i$
  with $\|\nabla f_i(\pmb x^*)\|_2 = \frac{1}{2}$.
\end{itemize}

\section*{C\quad Projection Onto Convex/Nonconvex Sets}

We list a few examples where the projection has closed-form expressions
or can be easily computed. These examples are frequently encountered
in science and engineering applications.
\begin{itemize}
\item The projection onto the linear subspace spanned by the columns of
the matrix $\pmb A\in\mathbb{R}^{n\times n_c}$ with $n_c\le n$ is
\begin{equation}\label{Proj:A}
    \Pi_{\pmb A}(\pmb z) = \pmb A\pmb A^\dagger\pmb z
\end{equation}
with $\pmb A^\dagger=(\pmb A^\top\pmb A)^{-1}\pmb A^\top$
being the Moore-Penrose pseudoinverse of $\pmb A$.

\item The projection onto hyperplane
$\mathcal{H}_i=\left\{\pmb x|\pmb a_i^\top\pmb x=b_i\right\}$ is
\begin{equation}\label{Proj:Hi}
    \Pi_{\mathcal{H}_i}(\pmb z)
    = \pmb z - \frac{\pmb a_i^\top\pmb z-b_i}{\|\pmb a_i\|_2^2}\pmb a_i
\end{equation}
with running time being $\mathcal{O}(n)$. The Kaczmarz method uses the projection
of \eqref{Proj:Hi} to update the iterate:
\begin{equation*}
    \pmb x^{k+1}= \pmb x^k - \frac{\pmb a_{i_k}^\top\pmb x^k-b_{i_k}}{\|\pmb a_{i_k}\|_2^2}\pmb a_{i_k}.
\end{equation*}
\item The projection onto circle $\mathcal{C}_i=\{\pmb x|\|\pmb x - \pmb c_i\|_2^2 = r_i^2\}$ is
\begin{equation}\label{Proj:Ci}
\Pi_{\mathcal{C}_i}(\pmb z) =
\left\{ \begin{array}{rc}
\pmb c_i + \frac{r_i}{\|\pmb z - \pmb c_i\|_2}(\pmb z - \pmb c_i),~\mathrm{if}~\pmb z \ne \pmb c_i \\
\pmb c_i + r_i \pmb v,~\mathrm{if}~\pmb z = \pmb c_i  \\
\end{array} \right.
\end{equation}
where $\pmb v$ is an arbitrary vector with unit norm $\|\pmb v\|_2=1$.
Computing $\Pi_{\mathcal{C}_i}(\pmb z)$ requires an $\mathcal{O}(n)$ running time.
\item The projection onto the surface of the ellipsoid
$\mathcal{E}_i=\left\{\pmb x\left||\pmb a_i^\mathrm{H}\pmb x|^2=b_i^2\right.\right\}$ is
\begin{equation}\label{Proj:Ei}
\Pi_{\mathcal{E}_i}(\pmb z) =
\left\{
\begin{array}{rc}
    \pmb z - \left(1-\frac{b_i}{|\pmb a_i^\mathrm{H}\pmb z|}\right)
    \frac{\pmb a_i^\mathrm{H}\pmb z}{\|\pmb a_i\|_2^2}\pmb a_i,
    ~\mathrm{if}~\pmb a_i^\mathrm{H}\pmb z \ne 0\\
    \pmb z - \frac{b_i}{\|\pmb a_i\|_2^2}\mathrm{e}^{\jmath\theta}\pmb a_i,
    ~\mathrm{if}~\pmb a_i^\mathrm{H}\pmb z = 0\\
\end{array} \right.
\end{equation}
where $\jmath$ is the imaginary unit and $\theta\in[0,2\pi)$ is an arbitrary phase angle.
Computing $\Pi_{\mathcal{E}_i}(\pmb z)$ takes time $\mathcal{O}(n)$.
\end{itemize}
The projections onto linear subspace and hyperplane are convex while those
onto the surfaces of the circle or ellipsoid are nonconvex. It is known that
the POCS is unique and non-expansive. The non-expansiveness of POCS refers to
\begin{equation}\label{Nonexpan:Convex}
    \|\Pi(\pmb z_1)-\Pi(\pmb z_2)\|_2\le\|\pmb z_1-\pmb z_2\|_2
\end{equation}
for any $\pmb z_1$ and $\pmb z_2$. However, the two properties do not necessarily
hold for projection onto nonconvex sets. For example, when $\pmb z = \pmb c_i$, the
projection onto the circle is not unique. It is also not non-expansive. When two
points are inside of the circle, the projection $\Pi_{\mathcal{C}_i}(\pmb z)$
of \eqref{Proj:Ci} can be expansive. That is, we have
$\|\Pi(\pmb z_1)-\Pi(\pmb z_2)\|_2>\|\pmb z_1-\pmb z_2\|_2$ for
$\|\pmb z_1 - \pmb c_i\|_2 < r_i$ and $\|\pmb z_2 - \pmb c_i\|_2 < r_i$
with $\pmb z_1\ne\pmb z_2$ and $(\pmb z_1 - \pmb c_i) \nparallel (\pmb z_2 - \pmb c_i)$.

\end{document}